\documentclass[11pt,reqno]{amsart} 
\usepackage{a4wide}
\usepackage{graphicx,amssymb,amsfonts,euscript,amsmath}

\newtheorem{m-theorem}[stuff]{Theorem}
\newtheorem{m-proposition}[stuff]{Proposition}
\newtheorem{m-corollary}[stuff]{Corollary}
\newtheorem{m-lemma}[stuff]{Lemma}
\newtheorem{m-definition}[stuff]{Definition}
\newtheorem{m-notation}[stuff]{Notation}
\newtheorem{m-remark}[stuff]{Remark}
\newenvironment{thm-nono}[1]{%theorem no-number
\vskip5pt\trivlist \itemindent 0pt
\item[\hskip\labelsep\bf Theorem #1.]%
\it\ignorespaces}{\endtrivlist\vskip5pt}% 
\newenvironment{m-proof}{%
\vskip4pt\trivlist \itemindent 0pt
\item[\hskip\labelsep\it Proof.]%
\ignorespaces}{\hfill$\Box$\endtrivlist\vskip4pt}%

\newcommand{\euf}{\EuScript}
\newcommand{\mfrak}{\mathfrak}
\newcommand{\eA}{{\euf A}}
\newcommand{\eF}{{\euf F}}
\newcommand{\eI}{{\euf I}}
\newcommand{\eL}{{\euf L}}
\newcommand{\eO}{{\euf O}}
\let\dashto\dashrightarrow
\let\mbb\mathbb
\newcommand{\Bl}{\mathop{\rm Bl}\nolimits}
\newcommand{\codim}{\mathop{\rm codim}\nolimits}
\let\ges\geqslant
\newcommand{\kk}{{\Bbbk}}
\let\les\leqslant
\let\nit\noindent
\newcommand{\Pic}{\mathop{\rm Pic}\nolimits}
\newcommand{\bloc}{\mathop{\mfrak{b}}} 
\newcommand{\sbloc}{\mathop{\rm{sb}}}
\newcommand{\surj}{\twoheadrightarrow}
\newcommand\uset[2]{{\displaystyle\mathop{\mbox{$#2$}}_{#1}}}

\let\unbar\underbar

\title[$q$-ampleness of the tensor product] %
{On the $q$-ampleness of the tensor product\\ of two line bundles}
\author{Mihai Halic, Roshan Tajarod}

\begin{document}

\begin{abstract}
We prove that the tensor product of two line bundles, one being $q$-ample and the other with sufficiently low-dimensional base locus, is still $q$-ample.
\end{abstract}

\maketitle

\section*{The result}\label{intro}

The goal of this note is to prove the following property of the $q$-ample cone of a projective variety. 

\begin{thm-nono}{A}\label{thm:main}
Let $X$ be a normal, irreducible projective variety defined over an algebraically closed field of characteristic zero. Consider $\eA,\eL\in\Pic(X)$ and denote the stable base locus of $\eA$ by $\sbloc(\eA)$. We assume that $\eL$ is $q$-ample and  
$$
\begin{minipage}{.95\textwidth}
\null\hfill 
$q\ges\dim\big(\sbloc(\eA)\big).$
\hfill $(\star)$
\end{minipage}
$$ 
Then $\eA\otimes \eL$ is $q$-ample too. 
\end{thm-nono}

The classes of the $q$-ample line bundles form an open cone in the vector space 
$$N^1(X)_{\mbb R}:=(\Pic(X)/\sim_{\rm num})\otimes_{\mbb Z}{\mbb R},$$ 
generated by invertible sheaves (line bundles) on $X$ modulo numerical equivalence (cf. \cite{dps,gk}). The tensor product of two $q$-ample line bundles is not $q$-ample in general (cf. \cite[Theorem 8.3]{totaro}), and therefore the $q$-ample cone is, usually, not convex. 

This situation contrasts the classical case of ample line bundles, corresponding to $q=0$, which generate a convex cone. Actually, it is well-known that the ample cone of a projective variety is stable under the addition of a numerically effective (nef) term.

Moreover, Sommese proved in \cite[Corollary 1.10.2]{so} that the tensor product of two \emph{globally generated}, $q$-ample line bundles is still $q$-ample. However, the concept of $q$-ampleness used in \textit{loc. cit.} is defined geometrically and it is based on the global generation of the line bundles. 

For this reason, it is natural to ask whether the $q$-ample cone is stable under the addition of suitable terms; by abuse of language, we call such a feature a `convexity property'. The theorem stated above can be viewed as an answer to this question.

% % % % % % % % % % % % % % % % % % % % % % % % % % % % % % % % % % % % % 

\section{Notation and proof}\label{sct:proof}

\begin{m-definition}\label{def:q}{\rm 
(cf. \cite[\S6]{totaro}) 
Let $X$ be a projective variety defined over an algebraically closed field $\kk$ of characteristic zero. A line bundle $\eL\in\Pic(X)$ is called \emph{$q$-ample} if, for all coherent sheaves $\eF$ on $X$, holds: 
\begin{equation}\label{eq:q}
\exists\,m_\eF\text{ such that }\forall\,m\ges m_\eF\;\forall\,t>q,
\quad H^t(X,\eF\otimes\eL^m)=0.
\end{equation} 
Since any coherent sheaf admits a finite resolution by locally free sheaves, it is enough to check the condition \eqref{eq:q} for $\eF$ locally free. 
}\end{m-definition}

The definition is closely related to the notions of $q$-positivity in \cite{dps} and of geometric $q$-ampleness in \cite{so}. These concepts are compared in \cite{mats}. 

Clearly, if $\eL$ is $q$-ample, then it is $q'$-ample, for all $q'\ges q$; the larger the value of $q$, the weaker the restriction on $\eL$. Any line bundle on $X$ is $\dim X$-ample; the first interesting case is $q=\dim X-1$. 
In \cite[Theorem 9.1]{totaro}, Totaro proved that the $(\dim X-1)$-ample cone is the complement in $N^1(X)_{\mbb R}$ of the negative of the closed effective cone. 

\begin{m-notation}{\rm 
For $\eA\in\Pic(X)$, we denote: 
\begin{enumerate}
\item 
$\bloc(\eA)$ the \emph{base locus} of $\eA$; it is the zero locus of the `universal section': 
$$
\begin{array}{l}
\eO_X\to H^0(X,\eA)^\vee\otimes\eA,\quad 
(x,1)\longmapsto\kern-2ex\uset{s\in\text{basis of $H^0(X,\eA)$}}{\sum}\kern-5ex %
s^\vee\otimes s(x),\;\;\forall\,x\in X, 
\\ 
\text{where ${(s^\vee)}_{s\in\text{basis of }H^0(X,\eA)}$ is the dual basis.}
\end{array}
$$
The scheme structure of $\bloc(\eA)$ is defined by the following sheaf of ideals: 
\begin{equation} \label{eq:I}
H^0(X,\eA)\otimes\eA^{-1}\surj\eI_{\bloc(\eA)}\subset\eO_X.
\end{equation}
The \emph{stable base locus} of $\eA$ is the closed subset of $X$ obtained as the set-theoretical intersection $\sbloc(\eA):=\uset{a\ges 1}{\bigcap}\bloc(\eA^a)_{\rm red};$ when $a$ is sufficiently large and divisible, $\sbloc(\eA)=\bloc(\eA^a)_{\rm red}$. 

\item 
$\kappa(\eA)$ the \emph{Kodaira-Iitaka dimension} of $\eA$; it is defined as:
$$
\begin{array}{rl}
\kappa(\eA)&:=\text{transcend.\,deg.}_\kk
\Big(\,\uset{a\ges 0}{\bigoplus} H^0(X,\eA^a)\Big)-1
\\ &\displaystyle\;
=\max_{a\ges1}\, 
\dim\big(\,{\rm Image}(X\dashto |\eA^a|)\,\big).
\end{array}
$$
The set $\{a\ges 1\mid H^0(X,\eA^a)\neq0\}$ is a semi-group under addition and consists, when $a$ is sufficiently large, of the multiples of a certain integer. Moreover, for $a\gg0$ in the set, the images of the rational maps $X\dashto |\eA^a|$ are birational to each other; in particular, their dimension is $\kappa(\eA)$. 
\end{enumerate}
For details, see \cite[Definition 2.1.3, Proposition 2.1.21, Theorem 2.1.33]{laz1}.
}\end{m-notation}

\begin{m-remark}{\rm
Related to our result, consider for instance the case $\dim(\sbloc(\eA))=0$, that is a power of $\eA$ is globally generated by its sections ($\eA$ is semi-ample); the $0$-ample cone is stable under the addition of a semi-ample term. 

At the other end of the scale, Totaro's result \cite[Theorem 9.1]{totaro} shows that the $(\dim X-1)$-ample cone is stable under the addition of line bundles $\eA=\eO_X(D)$, where $D$ is an effective divisor, that is those line bundles which admit a non-trivial section; clearly, in this case $\dim(\sbloc(\eA))\les\dim X-1$. 
}\end{m-remark}

\begin{m-lemma}\label{lm:kb}
Assume that the image of $X\dashto|\eA|$ is $\kappa(\eA)$-dimensional and  $\bloc(\eA)\neq\emptyset$. Then holds: $\;\kappa(\eA)\ges\codim_X\big(\bloc(\eA)\big)-1$. 
\end{m-lemma}

\begin{m-proof}
The equation \eqref{eq:I} implies that the blow-up $\Bl_\eI(X)$ of the ideal $\eI:=\eI_{\bloc(\eA)}$ is a closed subscheme the product $X\times\mbb P(H^0(X,\eA)^\vee)$; let $(\sigma,f)$ denote the inclusion morphism. For general $x\in\bloc(\eA)$, $\sigma^{-1}(x)$ is at least $\big(\codim_X\big(\bloc(\eA)\big)-1\big)$-dimensional and it is also contained in $f\big(\Bl_\eI(X)\big)$, which is $\kappa(\eA)$-dimensional. 
\end{m-proof}

Now we start proving the theorem A. 

\begin{m-proof} 
We observe that the statement is invariant after replacing $\eA$ by some power $\eA^a$: indeed, $(\eA\otimes\eL)^a=\eA^a\otimes\eL^a$ and the inequality $(\star)$ is preserved. Thus, henceforth, we may assume the following: 
$$
\sbloc(\eA)=\bloc(\eA)_{\rm red},\quad 
{\rm Image}\big(X\dashto|\eA|\,\big)\text{ is }\kappa(\eA)\text{-dimensional}. 
$$

If $\kappa(\eA)\ges1$, Bertini's theorem (cf. \cite{diaz+harb} \cite[Th\'eor\`eme 6.3]{jou}) implies that we have the exact sequences: 
\begin{equation}\label{eq:X}
0\to\eA^{-1}\otimes\eO_{X_{l-1}}\to\eO_{X_{l-1}}\to\eO_{X_l}\to0,
\quad l=1,\dots,\kappa(\eA),
\end{equation}
where 
\begin{equation}\label{eq:X2}
%\\[1ex]
\begin{array}{l}
X=:X_0\supset X_1\supset\dots\supset X_{\kappa(\eA)},\quad \dim X_l=\dim X-l, 
\\[1ex] 
X_l\in\big|{\rm Image}\big(H^0(\eA)\to H^0(\eA\otimes\eO_{X_{l-1}})\big)\big|
\text{ are very general,}
\\[1ex]
\kappa\big(\eA\otimes\eO_{X_{\kappa(\eA)}}\big)=0.
\end{array}
\end{equation} 
We distinguish two cases, whether $\sbloc(\eA)$ is empty or not.  
\medskip

\nit\unbar{\textit{Case $\sbloc(\eA)=\emptyset$}}\quad 
In this case, $\eA$ is globally generated and the image of $X\to|\eA|$ is $\kappa(\eA)$-dimensional. For shorthand, we denote $\kappa:=\kappa(\eA)$. 

We argue by descending induction on $q$. We will prove the following stronger statement: for all locally free sheaves $\eF$ on $X$, 
\begin{equation}\label{eq:ind-hyp}
\exists\,m_\eF\;\forall\,m\ges m_\eF\;\forall\,0\les j\les m\;\forall\,t>q,\quad 
H^t(X,\eF\otimes\eA^j\otimes\eF^m)=0.
\end{equation}
We fix such an $\eF$. If $\kappa=0$, then $\eA\cong\eO_X$ because it is globally generated, and there is nothing to prove. Thus we may assume $\kappa\ges 1$. 

Let $q=\dim X-1$. We tensor the exact sequence \eqref{eq:X}, with $l=1$, by $\eF\otimes\eA^j\otimes\eL^m$ and obtain, for $j=1,\dots,m$: 
$$
H^{\dim X}(X,\eF\otimes\eA^{j-1}\otimes\eL^m)\to
H^{\dim X}(X,\eF\otimes\eA^{j}\otimes\eL^m)\to0.
$$
Thus $H^{\dim X}(X,\eF\otimes\eL^m)\to H^{\dim X}(X,\eF\otimes\eA^j\otimes\eL^m)$ is surjective, and \eqref{eq:ind-hyp} follows. 

Suppose now that \eqref{eq:ind-hyp} holds for $q$, for some $m^{(q)}_\eF$, and let us prove it for $q-1$. So, if $\eL$ is $(q-1)$-ample (so it is $q$-ample), we must show that the $H^q(\cdot)$-term vanishes. The definition of the $(q-1)$-ampleness implies that there is $m^{(q-1)}_\eF\ges m^{(q)}_\eF$ such that 
\begin{equation}\label{eq:0}
H^q(X,\eF\otimes\eL^m\otimes\eO_{X_l})=0,\;\forall\,l=0,\dots,\kappa,\;\forall\,m\ges m^{(q-1)}_\eF.
\end{equation}
We observe that $\eA\otimes\eO_{X_\kappa}\cong\eO_{X_\kappa}$, because $\kappa(\eA\otimes\eO_{X_\kappa})=0$ and $\eA\otimes\eO_{X_\kappa}$ is globally generated, which implies: 
\begin{equation}\label{eq:kpa}
H^q(X,\eF\otimes\eA^j\otimes\eL^m\otimes\eO_{X_\kappa})=0,\text{ for }0\les j \les m. 
\end{equation}
Now assume that $X_l$ satisfies $H^q(X,\eF\otimes\eA^j\otimes\eL^m\otimes\eO_{X_l})=0$, for $0\les j \les m$,  and we prove the same for $X_{l-1}$. The sequence \eqref{eq:X} tensored by $\eF\otimes\eA^j\otimes\eL^m$, $j=1,\dots,m$, yields: 
$$
H^{q}(X,\eF\otimes\eA^{j-1}\otimes\eL^m\otimes\eO_{X_{l-1}})\to
H^{q}(X,\eF\otimes\eA^{j}\otimes\eL^m\otimes\eO_{X_{l-1}})\to0,\;\forall\,m\ges m^{(q-1)}_\eF. 
$$
The equation \eqref{eq:0} implies  $H^q(X,\eF\otimes\eA^j\otimes\eL^m\otimes\eO_{X_{l-1}})=0$, for $0\les j\les m$. Recursively, after $\kappa$ steps, we find that \eqref{eq:ind-hyp} holds for $X=X_0$ and $t=q$. This completes the inductive argument.
\medskip

\nit\unbar{\textit{Case $\sbloc(\eA)\neq\emptyset$}}\quad 
The key is again the exact sequences \eqref{eq:X}. Let 
$$\kappa:=\codim\big(\bloc(\eA)\big)-1.$$ 
The lemma \ref{lm:kb} implies that we have the inequality: $\;\kappa(\eA)\ges\kappa\ges\dim X-q-1.$ 
The term $X_\kappa$ in \eqref{eq:X}, has the following properties: 
\begin{itemize} %[itemindent=*]
\item 
$\kappa(\eA\otimes\eO_{X_\kappa})\ges0$; 
\item 
$\dim X_\kappa=\dim X-\kappa=\dim\big(\bloc(\eA)\big)+1$,\quad $\bloc(\eA)_{\rm red}\subset (X_\kappa)_{\rm red}$. 
\end{itemize}
Since the base locus is non-empty, there is a section in $\eA$ which vanishes along a (non-trivial) divisor $X_{\kappa+1}\subset X_\kappa$. (Otherwise, a component of $X_\kappa$ must be contained in $\bloc(\eA)$.) This yields one more exact sequence: 
$$
0\to\eA^{-1}\otimes\eO_{X_\kappa}\to\eO_{X_\kappa}\to\eO_{X_{\kappa+1}}\to0, 
\quad\dim X_{\kappa+1}=\dim X-(\kappa+1)\les q.
$$
We tensor it by $\eF\otimes\eA^j\otimes\eL^m$ and deduce, for $t>q$: 
$$
H^{t}(X,\eF\otimes\eA^{j-1}\otimes\eL^m\otimes\eO_{X_\kappa})\to
H^{t}(X,\eF\otimes\eA^{j}\otimes\eL^m\otimes\eO_{X_\kappa})\to0.
$$
Since $\eL$ is $q$-ample, it follows that for $t>q$ holds: 
$$
H^{t}(X,\eF\otimes\eA^{j}\otimes\eL^m\otimes\eO_{X_\kappa})=0,\quad
\text{ for }m\gg0,\;0\les j\les m.
$$
This is the vanishing \eqref{eq:kpa}, necessary for the induction step. Hence we can repeat the proof of the previous case.
\end{m-proof}

\begin{m-remark}{\rm 
\begin{enumerate}
\item 
The proof of Sommese's result \cite[Corrollary 1.10.2]{so} about the convexity of the cone generated by the geometrically $q$-ample line bundles does not carry over to our setting because it essentially uses their global generation. 
\item 
It is not clear to us whether the theorem remains valid if, instead of $q$-ample line bundles, one considers $q$-positive line bundles (cf. \cite{dps}). 
\end{enumerate}
}\end{m-remark}

\end{document}